\documentclass[number,secthm,citesort,seceqn,dvips]{arxbj}
\usepackage{multirow}
\usepackage{graphicx}

\aid{0}
\volume{17}
\issue{4}
\pubyear{2011}
\firstpage{1386}
\lastpage{1399}
\doi{10.3150/10-BEJ317}

\makeatletter

\newproclaim{definition}{Definition}[section]

\newremark{rem}{Remark}[section]

\newproclaim{assu}{Assumptions}
\makeatother

\begin{document}
\begin{frontmatter}

\title{Tensor-based projection depth}
\runtitle{Tensor-based projection depth}

\begin{aug}
\author{\fnms{Yonggang} \snm{Hu}\thanksref{e1}\ead[label=e1,mark]{xd7688@hotmail.com}\corref{}},
\author{\fnms{Yong} \snm{Wang}\thanksref{e2}\ead[label=e2,mark]{yongwang82@gmail.com}}
\and
\author{\fnms{Yi}
\snm{Wu}\thanksref{e3}\ead[label=e3,mark]{wuyi\_work@sina.com}}

\runauthor{Y.G. Hu, Y. Wang and Y. Wu}
\address{Department of Mathematics and Systems Science, National University of Defense Technology, Changsha 410073,
China.\\
\printead{e1,e2,e3}}
\end{aug}

\received{\smonth{5} \syear{2009}}
\revised{\smonth{6} \syear{2010}}

\begin{abstract}
The conventional definition of a depth function is vector-based. In
this paper, a novel projection depth (PD) technique directly based on
tensors, such as matrices, is instead proposed. Tensor projection depth
(TPD) is still an ideal depth function and its computation can be
achieved through the iteration of PD. Furthermore, we also discuss the
cases for sparse samples and higher order tensors. Experimental results
in data classification with the two projection depths show that TPD
performs much better than PD for  data with a natural tensor form, and
even when the data have a natural vector form, TPD appears to perform
no worse than PD.
\end{abstract}


\end{frontmatter}

\section{Introduction}\label{sec:introduction}

In the last ten years, statistical depth functions have  increasingly
served as a useful tool in multidimensional exploratory data analysis
and inference. The depth of a point in the multidimensional space
measures the centrality of that point with respect to a multivariate
distribution or a given multivariate data cloud. Depth functions have
been successfully used in many fields, such as quality indices
\cite{Liu1992,Liu1993}, multivariable regression \cite{Tian2002},
limiting $p$ values \cite{Liu1997}, robust estimation \cite{Chen2004},
nonparametric tests \cite{Chenouri2004} and discriminant analysis
\cite{Ghosh2005a,Ghosh2005b,Jornsten2004,Cui2008}. Some common
statistical depths which have been defined include  \textit{half-space
depth} \cite{Tukey1975}, \textit{simplicial depth} \cite{Liu1990},
\textit{projection depth}
\cite{Stahel1981,Donoho1982,Donoho1992,Zuo2003}, \textit{spatial depth}
\cite{Vardi2000}, \textit{spatial rank depth} \cite{Gao2003} and
\textit{integrated dual depth} \cite{Cuevas2009}. Compared to the
others, projection depth (PD) is preferable because of its good
properties such as robustness, affine invariance, maximality at center,
monotonicity relative to deepest point, vanishing at infinity and so
on.

However, almost all the depths proposed in the literature are defined
over the vector space by now, and the fact is that not all of the
observations are naturally in vector form. In the real world, the
extracted feature of an object often has some specialized structures,
and such structures are in the form of a second, or even higher order
tensor. For example, this is the case when a captured image is a
second-order tensor, that is, a matrix, and when the sequential data,
such as a video sequence for event analysis, is in the form of a
third-order tensor. It would be desirable to keep the underlying
structures of the data unchanged  during the data analysis.

Most of the previous work on depth has first transformed the input
tensor data into vectors, which in fact changes the underlying
structure of the data sets. At the same time, such a transformation
often leads to the curse of dimensionality problem and the small sample
size problem since most depth functions (such as Mahalanobis depth)
require the covariance matrix to be positive definite.

Therefore, it is necessary to extend the definition of depth to tensor
spaces in order to process the data sets directly with tensors without modifying
the structures of them. In fact, many tensor-based methods in
discriminant analysis have been proposed and have led to many nice
results \cite{Cai2006,Xu2006,Zhi2008}. In this paper, informed by the
aforementioned works, we propose a tensor-based projection depth (TPD)
in order to extend the definition of projection depth to  tensor
spaces. We will prove that TPD is still an ideal depth according to the
criteria \cite{Zuo2000a}. Also, we will explore the characteristics of
high order tensor projection depth in theory. We will demonstrate that
TPD allows us to avoid the above two problems when using vector
representation.

The paper is organized as follows. Section \ref{sec:tensor algebra}
briefly introduces tensor algebra. Section \ref{sec:projection depth}
introduces the projection depth and gives the solution to the Rayleigh
projection depth. Section \ref{sec:tensor projection depth} gives the
definition of tensor projection depth and discusses its properties.
Section \ref{sec:algorithm} supplies the algorithm for TPD and analyzes
its convergence. Section \ref{sec:for sparse samples}  analyzes the
special case of sparse samples. Section \ref{sec:from matrix to higher
tensor}  discusses the TPD for higher order tensors. Section
\ref{sec:experiments} gives  numerical results for TPD. Section
\ref{sec:conclusion} concludes the paper, and  proofs of selected
theorems and propositions are given in the Appendix.

\section{Tensor algebra}
\label{sec:tensor algebra}

A tensor $T$ of order $k$ is a real-valued multilinear function on $k$
vector spaces \cite{Itskov2007}:
\[
  T\dvtx \mathbb{R}^{n_1}\times\cdots\times\mathbb{R}^{n_k}\rightarrow\mathbb{R.}
\]
A multilinear function is linear as a function of each variable
considered separately. The set of all $k$th-order tensors on
$\mathbb{R}^{n_i}$, $i=1, \ldots,k$, denoted by $\mathcal{T}^k$, is a
vector space under the usual operations of pointwise addition and
scalar multiplication:
\begin{eqnarray*}
  (aT)(\mathbf{a}_1,\ldots,\mathbf{a}_k)&=&a(T(\mathbf{a}_1,\ldots,\mathbf{a}_k)),\\
  (T+T')(\mathbf{a}_1,\ldots,\mathbf{a}_k)&=&T(\mathbf{a}_1,\ldots,\mathbf{a}_k)+ T'(\mathbf{a}_1,\ldots,\mathbf{a}_k),
\end{eqnarray*}
where $\mathbf{a}_i\in\mathbb{R}^{n_i}$.

Given two tensors, $S\in \mathcal{T}^k$ and $T\in \mathcal{T}^l$, their
product,
\[
  S\otimes T\dvtx \mathbb{R}^{n_1}\times\cdots\times\mathbb{R}^{n_{k+l}} \rightarrow\mathbb{R},
\]
is defined as
\[
  S\otimes T(\mathbf{a}_1,\ldots,\mathbf{a}_{k+l})= S(\mathbf{a}_1,\ldots,\mathbf{a}_k)T(\mathbf{a}_{k+1},\ldots,\mathbf{a}_{k+l}).
\]
It is immediate from the multilinearity of $S$ and $T$ that $S\otimes
T$ depends linearly on each argument $\mathbf{a}_i$ separately, so it
is a ($k+l$)th-order tensor.

First-order tensors are simply  vectors on $\mathbb{R}^{n_1}$. That is,
$\mathcal{T}_1 =\mathcal{R}^{n_1}$, where $\mathcal{R}^{n_1}$ is the
dual space of $\mathbb{R}^{n_1}$. A second-order tensor space is a
product of two first-order tensor spaces, that is, $\mathcal{T}^2 =
\mathcal{R}^{n_1}\otimes \mathcal{R}^{n_2} $. Let
$\mathbf{e}_1,\ldots,\mathbf{e}_{n_1}$ be the standard basis of
$\mathbb{R}^{n_1}$ and $\varepsilon_1,\ldots,\varepsilon_{n_1}$ be the
dual basis \cite{Meise1997} of $\mathcal{R}^{n_1}$ which is formed from
coordinate functions with respect to the basis of $\mathcal{R}^{n_1}$.
Likewise, let $\tilde{\mathbf{e}}_1,\ldots,\tilde{\mathbf{e}}_{n_1}$ be
a basis of $\mathbb{R}^{n_2}$  and
$\tilde{\varepsilon}_1,\ldots,\tilde{\varepsilon}_{n_1}$ be the dual
basis of $\mathcal{R}^{n_2}$. We have
\[
  \varepsilon_i(\mathbf{e}_j)=\delta_{ij}\quad\mbox{and} \quad\tilde{\varepsilon}_i(\tilde{\mathbf{e}}_j)=\delta_{ij},
\]
where $\delta_{ij}$ is the Kronecker delta function. Thus,
$\{\varepsilon_i\otimes\tilde{\varepsilon}_j\}$ $(1\leq i\leq n_1,
1\leq j\leq n_2)$ forms a basis for
$\mathcal{R}^{n_1}\otimes\mathcal{R}^{n_2}$. For any second-order
tensor $T$, we can write
\[
  T=\sum_{i,j}
  T_{ij}\varepsilon_i\otimes\tilde{\varepsilon}_j.
\]

Given two vectors
$\mathbf{a}=\sum_{k=1}^{n_1}a_k\mathbf{e}_k\in\mathbb{R}^{n_1}$ and
$\mathbf{b}=\sum_{l=1}^{n_2}b_l\tilde{\mathbf{e}}_l\in\mathbb{R}^{n_2}$,
we have
\begin{eqnarray}\label{eq:tensor with matrix}
T(\mathbf{a},\mathbf{b})
&=&
\sum_{ij}T_{ij}\varepsilon_i\otimes\tilde{\varepsilon}_j\Biggl(\sum_{k=1}^{n_1}a_k\mathbf{e}_k,\sum_{l=1}^{n_2}b_l\tilde{\mathbf{e}}_l\Biggr)\nonumber
\\
&=&
\sum_{ij}T_{ij}\varepsilon_i\Biggl(\sum_{k=1}^{n_1}a_k\mathbf{e}_k\Biggr)\tilde{\varepsilon}_j\Biggl(\sum_{l=1}^{n_2}b_l\tilde{\mathbf{e}}_l\Biggr)
\\
&=&\sum_{ij}T_{ij}a_ib_j=\mathbf{a}^TT\mathbf{b}.\nonumber
\end{eqnarray}
This shows that every second-order tensor in
$\mathcal{R}^{n_1}\otimes\mathcal{R}^{n_2}$ uniquely corresponds to an
$n_1\times n_2$ matrix.

Note that in this paper, our primary interest is focused on
second-order tensors. However, most of our conclusions for
second-order TPD can be naturally extended to higher orders. We will
discuss this question in  Section \ref{sec:from matrix to higher
tensor}.

\section{Projection depth}\label{sec:projection depth}

According to  \cite{Zuo2003}, the definition of projection depth can be
expressed as follows.

\begin{definition}\label{def:projection depth}
Let $\mu$ and $\sigma$ be univariate
location and scale measures, respectively. Define the outlyingness of a
point $\mathbf{x}\in \mathbb{R}^p$ with respect to a given function $F$
of $X$ in $\mathbb{R}^p$, $p\geq 1$, as
\begin{equation}\label{eq:outlying_vector}
O(\mathbf{x},F)=\sup_{\|\mathbf{u}\|=1}
\frac{|\mathbf{u}^T\mathbf{x}-\mu(F_\mathbf{u})|}{\sigma(F_{\mathbf{u}})},
\end{equation}
where $F_\mathbf{u}$ is the distribution of $\mathbf{u}^TX$. Then,
$O(\mathbf{x},F)$ is defined to be $0$ if
$\mathbf{u}^T\mathbf{x}-\mu(F_\mathbf{u})=\sigma(F_{\mathbf{u}})=0$.
The projection depth (PD) of a point $\mathbf{x}\in \mathbb{R}^p$ with
respect to the given $F, \mathit{PD}(\mathbf{x},F)$, is then defined as
\begin{equation}
\mathit{PD}(\mathbf{x}, F)=\frac{1}{1+O(\mathbf{x}, F)}.
\end{equation}
\end{definition}

\begin{rem}
Here, we also assume that $\mu$ and $\sigma$ exist uniquely, $\mu$ is
translation  and scale equivariant, and $\sigma$ is scale equivariant
and translation invariant, that is, $\mu(F_{sY+c})=s\mu(F_Y)+c$ and
$\sigma(F_{sY+c})=|s|\sigma(F_Y)$, respectively, for any scalars $s$,
$c$ and random variable $Y\in \mathbb{R}^1$.
\end{rem}

The most popular outlying function is defined as
\begin{equation}\label{eq:projection core}
O(\mathbf{x},F)=\sup_{\|\mathbf{u}\|=1}\frac{|\mathbf{u}^T\mathbf{x}-\mathit{Med}(F_\mathbf{u})|}
{\mathit{MAD}(F_\mathbf{u})},
\end{equation}
where $F_\mathbf{u}$ is the distribution of $\mathbf{u}^TX,
\mathit{Med}(F_\mathbf{u})$ is the median of $F_{\mathbf{u}}$ and
$\mathit{MAD}(F_\mathbf{u})$ is the median of the distribution of
$|\mathbf{u}^TX-\mathit{Med}(F_u)|$.

Apart from the good properties of a statistical depth function, this
version of PD is more robust compared with other depths. However, it is
hard to compute for  high-dimensional samples.

Obviously,  the variance and mean are also natural choices for $\sigma$
and $\mu$, respectively. It is easy to prove that such a
projection-based depth is also an ideal depth function. And, most
importantly, its computation is very simple.

\begin{thm}[(Rayleigh projection depth)]
\label{thm:rayleigh solution} Let $(\mu,\sigma)={}$(\textup{mean},
\textup{variance}), and suppose that the second moments of $X$ exist
and that $X\sim F$. The solution of the outlying function
(\ref{eq:outlying_vector}) is then that of  a Rayleigh quotient
problem,
\begin{eqnarray}\label{eq:rayleigh_projection depth}
O_R(\mathbf{x},F)
&=&
\sup_{\|\mathbf{u}\|=1}\frac{|\mathbf{u}^T\mathbf{x}-E(\mathbf{u}^TX)|} {\sqrt{E(\mathbf{u}^TX-E(\mathbf{u}^TX))^2}}\nonumber
\\[-10pt]\\[-10pt]
&=&
\sqrt{\frac{\mathbf{u}_1^T\mathbf{A}\mathbf{u}_1}{\mathbf{u}_1^T\mathbf{B}\mathbf{u}_1}}=\sqrt{\lambda_1},\nonumber
\end{eqnarray}
where $\mathbf{A}$ is the matrix $(\mathbf{x}-EX)(\mathbf{x}-EX)^T$,
$\mathbf{B}$ is $E(X-EX)(X-EX)^T$, $\lambda_1$ is the largest
eigenvalue of the generalized eigenvalue problem
\[
\mathbf{A}\mathbf{z}=\lambda\mathbf{B}\mathbf{z},\qquad\mathbf{z}\neq
0,
\]
and $\mathbf{u}_1$ is the corresponding eigenvector of
$\lambda_1$.

We call this projection depth the \textit{Rayleigh projection depth}.
\end{thm}

\begin{rem}
In this paper, for the convenience of computation, the examples in the
experiments are all based on the Rayleigh projection depth, that is,
$(\mu,\sigma)={}$(\textup{mean}, \textup{varance}).
\end{rem}

\begin{rem}
Obviously, RPD requires the covariance $\mathbf{B}$ to be positive.
%
To avoid this situation, for the sparse samples, we simply project the
samples into their nonzero subspace using principal component analysis
(PCA).
\end{rem}

\section{Tensor projection depth}\label{sec:tensor projection depth}

Before describing tensor projection depth, we first review the
terminology associated with tensor operations
\cite{Kolda2001,Lathauwer2000}. The inner product of tensors
$\mathbf{A}$ and $\mathbf{B}$ (with the same orders and dimensions) is
$\langle \mathbf{A}, \mathbf{B}\rangle =
\sum_{i,j}\mathbf{A}_{ij}\mathbf{B}_{ij}$. The norm of a tensor
$\mathbf{A}$ is defined as its Frobenius norm, that is,
$\|\mathbf{A}\|=\sqrt{\langle \mathbf{A}, \mathbf{A}\rangle}$, and the
distance between two tensors $\mathbf{A}$ and $\mathbf{B}$ in
$\mathcal{R}^{n_1}\otimes\mathcal{R}^{n_2}$ is defined as
$\|\mathbf{A}-\mathbf{B}\|$, where
$\mathbf{A}-\mathbf{B}=(\mathbf{A}_{ij}-\mathbf{B}_{ij})_{n_1\times
n_2}$.

From the tensorial viewpoint, if we take $X$ as a random variable in
the first-order tensor space $\mathcal{R}^{n_1}$, then the outlyingness
of the projection depth in Definition \ref{def:projection depth} can be
expressed as
\[
O(\mathbf{x},X)=\sup_{\|\mathbf{u}\|=1}\frac{|\mathbf{x}(\mathbf{u})-\mu(X(\mathbf{u}))|}{\sigma(X(\mathbf{u}))}.
\]

Thus, if $\mathcal{X}\in \mathcal{R}^{n_1}\otimes\mathcal{R}^{n_2}$ is
a random variable, then, according to the formula (\ref{eq:tensor with
matrix}), the outlying function in the tensor space
$\mathcal{R}^{n_1}\otimes\mathcal{R}^{n_2}$ can be naturally defined as
\begin{equation}
\label{eq:projection gen core tensor}
O(\mathbf{X},\mathcal{X})=\sup_{\|\mathbf{u}\|=\|\mathbf{v}\|=1}
\frac{|\mathbf{X}(\mathbf{u},\mathbf{v})-\mu(\mathcal{X}(\mathbf{u,v}))|}
{\sigma(\mathcal{X}(\mathbf{u,v}))}=\sup_{\|\mathbf{u}\|=\|\mathbf{v}\|=1}
\frac{|\mathbf{u}^T\mathbf{X}\mathbf{v}-\mu(\mathbf{u}^T\mathcal{X}\mathbf{v})|}
{\sigma(\mathbf{u}^T\mathcal{X}\mathbf{v})},
\end{equation}
where $\mathbf{u}\in \mathbb{R}^{n_1}$ and $\mathbf{v}\in
\mathbb{R}^{n_2}$.

\begin{definition}[(Tensor projection depth)]
The projection depth with  outlying function given by formula
(\ref{eq:projection gen core tensor}) is called \textit{tensor
projection depth}.
\end{definition}

For a given univariate location (or ``center'') measure $\mu$, a
distribution function $F_\mathcal{X}$ is called $\mu$-\textit{symmetric
about the point} $\theta \in \mathcal{R}^{n_1}\otimes\mathcal{R}^{n_2}$
if
$\mu(\mathbf{u}^T\mathcal{X}\mathbf{v})=\mathbf{u}^T\theta\mathbf{v}$
for any pair of unit vectors $\mathbf{u}\in \mathbb{R}^{n_1},
\mathbf{v}\in \mathbb{R}^{n_2}$. We have the following theorem.

\begin{thm}
\label{thm:properties}
Suppose that $\theta$ in
$\mathcal{R}^{n_1}\otimes\mathcal{R}^{n_2}$ is the point of symmetry of
a distribution $F(\mathcal{X})$ with respect to a given notion of
symmetry. The tensor projection depth function
$\mathit{TPD}(\mathbf{X},\mathcal{X})$ is:
\begin{enumerate}
\item convex;
\item symmetric for $\mu$-symmetric $F$;
\eject
\item affine invariant;
\item monotonic relative to the deepest point;
\item vanishing at infinity, that is, $\mathit{TPD}(\mathbf{X},\mathcal{X})\rightarrow 0$ as $\|\mathbf{X}\|\rightarrow \infty$;
\item maximized at the center of $\mu$-symmetric $F$.
\end{enumerate}
\end{thm}

\begin{rem}
Theorem \ref{thm:properties} shows that TPD is still an ideal depth
according to the criteria \cite{Zuo2000a}. Furthermore, we can easily
obtain many other properties of TPD beyond those of the PD in
\cite{Zuo2003}, such as the properties of its sample versions and its
medians. However, these are not the key points of this paper and so we
omit any detailed discussion here.
\end{rem}

\section{Algorithm}
\label{sec:algorithm}
Suppose that the elements of
$S_n=\{\mathbf{X}_1,\dots,\mathbf{X}_n\}$ are generated from $F$ (where
$F_n$ is its empirical distribution) and that $\mathbf{X}$ is a fixed
tensor. The TPD of $\mathbf{X}$ with respect to $F_n$ can then be
computed by the following algorithm:
\begin{enumerate}
\item \textit{Initialization:} Let $\mathbf{u}=(1,\ldots,1)^T$.
\item\textit{Computing $\mathbf{v}$}: Let
$\mathbf{x}_i=\mathbf{X}_i^T\mathbf{u}$ and
$F^{\mathbf{u}}_n=\mathbf{u}^TF_n$. Then, $\mathbf{v}$ can be computed
by solving the vector-based projection depth
\begin{equation}
\label{eq:algorithm step 2}
\sup_{\|\mathbf{v}\|=1}\frac{|\mathbf{v}^T\mathbf{x}- \mu(F_n^{\mathbf{u}}\mathbf{v})|}{\sigma(F_n^{\mathbf{u}}\mathbf{v})}.
\end{equation}
\item \textit{Computing $\mathbf{u}$:} Once $\mathbf{v}$ is obtained,
let $\tilde{\mathbf{x}}_i=\mathbf{X}_i\mathbf{v}$ and
$F^{\mathbf{v}}_n=F_n\mathbf{v}$. Then, $\mathbf{u}$ can be computed by
solving the following optimization problem:
\begin{equation}
\label{eq:algorithm step 3}
\sup_{\|\mathbf{u}\|=1}\frac{|\mathbf{u}^T\tilde{\mathbf{x}}-
\mu(\mathbf{u}^TF^{\mathbf{v}}_n)|}
{\sigma(\mathbf{u}^TF^{\mathbf{v}}_n)}.
\end{equation}
\item \textit{Iteratively computing $\mathbf{u}$ and $\mathbf{v}$}:
Using steps 2 and 3, we can iteratively compute $\mathbf{u}$ and
$\mathbf{v}$ until they tend to converge.
\end{enumerate}

\begin{rem}
The optimization problems (\ref{eq:algorithm step 2}) and
(\ref{eq:algorithm step 3}) are the same as (\ref{eq:projection core})
in the vector-based projection depth algorithm. Thus, any computational
method for the projection depth can also be used here.
\end{rem}

The following theorem shows that the above algorithm converges.

\begin{thm}\label{thm:algorithm convergence}
The iterative procedure to solve the optimization problems
(\ref{eq:algorithm step 2}) and (\ref{eq:algorithm step 3}) will
monotonically increase the objective function value in
(\ref{eq:projection gen core tensor}),  hence the algorithm converges.
\end{thm}

\begin{rem}
Furthermore, if the optimization problem  (\ref{eq:outlying_vector}) is
convex, then the solution of (\ref{eq:projection gen core tensor}) is
also globally optimal. For instance, if $(\mu,\sigma)={}$(mean,
variance) (the Rayleigh projection depth), then its solution is also
globally optimal.
\end{rem}

\section{Sparse samples}
\label{sec:for sparse samples} As with RPD, TPD based on RPD also faces
the problem of sparse samples. From formulas (\ref{eq:algorithm step
2}) and (\ref{eq:algorithm step 3}), we know that for any sample set
$S_n=\{\mathbf{X}_1, \dots, \mathbf{X}_n\}$ and its corresponding
empirical distribution $F_n$, the algorithm in the previous section
requires the covariance matrices of $F^{\mathbf{u}}_n$ and
$F^{\mathbf{v}}_n$ to be positive for any $\mathbf{u}\in
\mathbb{R}^{n_1}$ and $\mathbf{v}\in \mathbb{R}^{n_2}$. However, in
practice, the tensor data usually do not satisfy such requirements.

There are two factors that can lead to such non-positiveness. First,
the sample size is too small, that is, the size of $S_n$ is less than
$n_1$ or $n_2$. Second, the data have some common columns or rows
(e.g., the images have identical color edges or patterns). In the
vector space, we usually use PCA to remove the redundant null space of
the samples and therefore we can use the tensor PCA proposed by Cai
\textit{et al.} \cite{Cai2005} to reduce the dimensionality of the
tensor samples.

Suppose that $M_\mathbf{X}=\frac{1}{n}\sum_{i=1}^n\mathbf{X}_i$,
\begin{eqnarray*}
 M_V&=&\sum_{i=1}^n\bigl((\mathbf{X}_i-M_{\mathbf{X}})(\mathbf{X}_i-M_{\mathbf{X}})^T\bigr),
  \\
 M_U&=&\sum_{i=1}^n\bigl((\mathbf{X}_i-M_{\mathbf{X}})^T(\mathbf{X}_i-M_{\mathbf{X}})\bigr),
\end{eqnarray*}
where the columns of $V$ are the eigenvectors of $M_V$, and $U$ are the
eigenvectors of $M_U$. Thus, the new mappings of $F_n$ can be expressed
as
\begin{equation}\label{eq:tensor_mappings}
F_n^{(r_1,r_2)}=\{V_{r_1}^T\mathbf{X}_1U_{r_2},\dots,V_{r_1}^T\mathbf{X}_nU_{r_2}\},
\end{equation}
where $r_1$ and $r_2$ are the mapping dimensions, and $V_{r_1}$ and
$U_{r_2}$ are the first $r_1$ and $r_2$ columns of $V$ and $U$,
respectively. Here, we take $r_1$ and $r_2$ to be the ranks of $M_V$
and $M_U$.



\begin{thm}
For any $\mathbf{u}\in \mathbb{R}^{r_1}$, $\mathbf{v}\in
\mathbb{R}^{r_2}$ with  $\|\mathbf{u}\|=\|\mathbf{v}\|=1$, the
covariance matrices of $\mathbf{u}^TF_n^{(r_1,r_2)}$ and
$F_n^{(r_1,r_2)}\mathbf{v}$ are always positive.
\end{thm}

\section{Higher order tensors}\label{sec:from matrix to higher tensor}

The algorithm described above takes second-order tensors (i.e.,
matrices) as input data. However, the algorithm can also be extended to
higher order tensors. In this section, we briefly describe the TPD
algorithm for  higher order tensors.

Let $S_n=\{\mathbf{X}_i, i=1,\ldots, n\}$ denote the sample set and
$F_n$  its empirical distribution, where $\mathbf{X}_i\in
\mathcal{R}^{n_1}\otimes\cdots\otimes\mathcal{R}^{n_k}$. The outlying
function of TPD is then
\begin{equation}\label{eq:core high order tensor}
O(\mathbf{X}, \mathcal{X})\doteq
\sup_{\|\mathbf{u}_1\|=\cdots=\|\mathbf{u}_k\|=1}
\frac{|\mathbf{X}(\mathbf{u}_1,\ldots,\mathbf{u}_k)-\mu(F_n(\mathbf{u}_1,
\ldots, \mathbf{u}_k))|} {\sigma(F_n(\mathbf{u}_1, \ldots,
\mathbf{u}_k))},
\end{equation}
where $\mathbf{u}_i\in \mathbb{R}^{n_i}$.

Before stating the algorithm, we first introduce an item of notation
which we will need. If $T\in
\mathcal{R}^{n_1}\otimes\cdots\otimes\mathcal{R}^{n_k}$, then for any
$\mathbf{a}_l\in \mathbb{R}^{n_l}$, $1\leq l\leq k$, we use
$T\times_{l} \mathbf{a}_{l}$ to denote a new tensor in
$\mathcal{R}^{n_1}\otimes\cdots\otimes\mathcal{R}^{n_{l-1}}\otimes\mathcal{R}^{n_{l+1}}\otimes\cdots\otimes\mathcal{R}^{n_k}$,
namely
\begin{equation}
T\times_{l}
\mathbf{a}_{l}=\sum_{i_l=1}^{n_l}T_{i_1,\ldots,i_{l-1},\ldots,i_{l+1},
\ldots,i_k}\cdot a_{i_l}.
\end{equation}
Thus, the algorithm for  higher order tensors can naturally be
expressed as follows:
\begin{enumerate}
\item \textit{Initialization}: Let
$\mathbf{u}_i^0=(x_1,\ldots,x_{n_i})^T$, $x_j\in\mathbb{R}$,
$j=1,\ldots,n_i$, $i=1,\ldots,k-1$.
\item \textit{Computing
$\mathbf{u}_k^0$}: If we let
$\mathbf{x}^k=\mathbf{X}\times_1\mathbf{u}_1^0\times_2\mathbf{u}_2^0\times\cdots\times_{k-1}\mathbf{u}_{k-1}^0$,
then $\mathbf{u}_k^0$ can be computed by solving the vector-based
projection depth
\begin{equation}\label{eq:algo highorder step 2}
\sup_{\|\mathbf{u}_k^0\|=1}\frac{|{\mathbf{u}_k^0}^T\mathbf{x}^k-
\mu(F_n\times_1\mathbf{u}_1^0\times_2\mathbf{u}_2^0\times\cdots\times_{k-1}\mathbf{u}_{k-1}^0)|}
{\sigma(F_n\times_1\mathbf{u}_1^0\times_2\mathbf{u}_2^0\times\cdots\times_{k-1}\mathbf{u}_{k-1}^0)}.
\end{equation}
\item \textit{Computing $\mathbf{u}_{k-1}^1$}: Once $\mathbf{u}_k^0$ is
obtained, we let
$\mathbf{x}^{k-1}=\mathbf{X}\times_1\mathbf{u}_1^0\times\cdots\times_{k-2}\mathbf{u}_{k-2}^0\times_k\mathbf{u}_k^0$
and $\mathbf{u}_{k-1}^1$ can be computed by solving the  optimization
problem
\begin{equation}
\label{eq:algo highorder step 3}
\sup_{\|\mathbf{u}_{k-1}^1\|=1}\frac{|{\mathbf{u}_{k-1}^{1T}}\mathbf{x}^{k-1}-
\mu(F_n\times_1\mathbf{u}_1^0\times\cdots\times_{k-2}\mathbf{u}_{k-2}^0\times_k\mathbf{u}_k^0)|}
{\sigma(F_n\times_1\mathbf{u}_1^0\times\cdots\times_{k-2}\mathbf{u}_{k-2}^0\times_k\mathbf{u}_k^0)}.
\end{equation}
\item \textit{Iteratively computing $\mathbf{u}_i$, $i=1,\ldots,k,$}
until they tend to converge.
\end{enumerate}

\begin{rem}
It is easy to prove that TPD in a higher order tensor space still
satisfies the above theorems and that its convergence is also
guaranteed by Theorem \ref{thm:algorithm convergence}.
\end{rem}

\section{Experiments}\label{sec:experiments}

First, we use data classification to
demonstrate the validity of the TPD. Consider a multivariate data set
$C$ that is partitioned into given classes $C_1,\dots,C_q$. An
additional data point $\mathbf{x}$ has to be assigned to one of several
given classes of \textit{object}. Suppose that there are $q$ classes.
The most natural classifier provided by \cite{Jornsten2004} is
then
\begin{equation} \label{eq:classifier}
\mathit{classd}(\mathbf{x})=\arg\max_j D(\mathbf{x}|C_j),
\end{equation}
where $D(\mathbf{x}|C_j)$ is the depth of the $\mathbf{x}$ with respect
to class $C_j$, $i=1,\dots,q$. This assigns $\mathbf{x}$ to the class
$C_j$ in which $\mathbf{x}$ is deepest.

The Columbia Object Image Library (COIL-$20$) \cite{Coil20} is a
database of grayscale images of $20$ objects. The objects were placed
on a motorized turntable against a black background. The turntable was
rotated through $360$ degrees to vary the object pose with respect to a
fixed camera. Images of the objects were taken at pose intervals of $5$
degrees. This corresponds to $72$ images with the dimensions of
$32\times 32$ pixels per object. Here, we only take the first $10$ objects as
examples.

\begin{figure}

\includegraphics{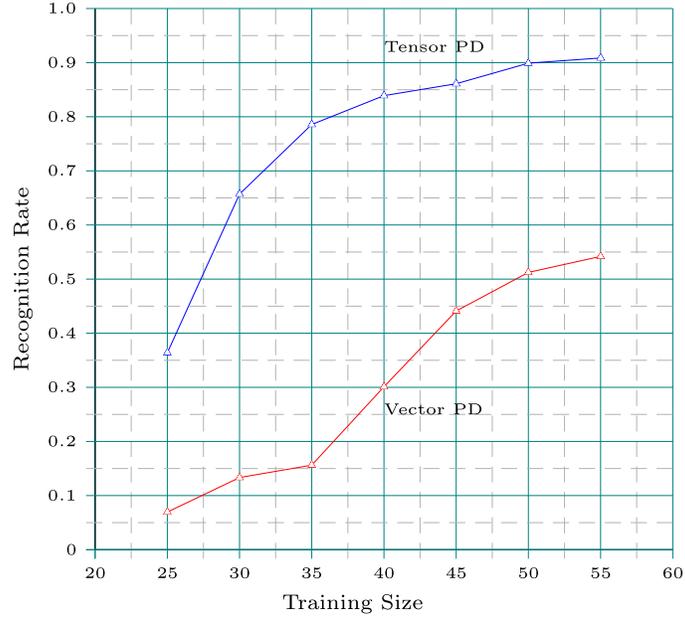}

\caption{Recognition rates of TPD and PD under different training
sample sizes.} \label{fig:reogntionrate_coil}
\vspace*{-12pt}
\end{figure}

In the experiments,  recognition rates under different training sizes
are computed by means of the following steps:
\begin{enumerate}
\item \textit{Select the test sets}: Randomly select $p$ test sets
$\mathbf{X}_{\mathrm{test}}^j$ from the object set $\mathbf{X}_j$ for
each class, where $j=1,\dots,10$.
\item \textit{for} each training size
$n_k$
\item[] \quad \textit{for} each repeating round $t$:
\begin{itemize}
\item \textit{Randomly select the training sets}: Randomly select $n_k$
training sets from $\mathbf{X}_j/\mathbf{X}_{\mathrm{test}}^j$ (the
left samples of $\mathbf{X}_j$) for each $j$, $j=1,\dots,10$. \item
\textit{Compute the recognition rate}. Compute the correctly recognized
number $\ell_j$ for each test set $\mathbf{X}_{\mathrm{test}}^j$ by
using the formula (\ref{eq:classifier}) and compute the glossary
recognition rate by $\eta_t=\sum_{j=1}^{10}\ell_j/10p$.
\end{itemize}
\item Compute the mean and variance of $\eta_t$.
\end{enumerate}
Here, $p=7$ and the training number equals 25, 30, 35, 40, 45, 50, 55,
respectively. The results are shown in Figure
\ref{fig:reogntionrate_coil} and Table \ref{tab:coil_PD_TPD}.

From Figure \ref{fig:reogntionrate_coil} and Table
\ref{tab:coil_PD_TPD}, we can see that for such samples with intrinsic
tensor form, TPD performs better than PD. A  question then naturally
arises: If the data sets are naturally in vector form, how does TPD
perform compared with PD? We will answer the question by means of the
following experiment.

We consider the famous Iris data \cite{fis}, which contains
measurements of four different features (sepal length, sepal width,
petal length and petal width) for each of 150 observations from three
different types of iris plant: (1) setosa; (2) virginica; (3)
versicolor. We randomly choose 10 observations from each class to
construct the test sets and then randomly select 10, 15, 20, 25, 30,
35, 40 samples from the remaining observations as the respective
training sets. For the computation of TPD, the samples are reshaped as
$2\times 2$.

\begin{table}
\caption{The mean, deviation and variance of the recognition rates by
TPD and PD with the COIL-20 set} \label{tab:coil_PD_TPD}
\begin{tabular*}{\textwidth}{@{\extracolsep{\fill}}lllllllll@{}}
\hline
\multirow{2}{31pt}{Training} & \multicolumn{4}{l}{Tensor projection depth}&\multicolumn{4}{l@{}}{Projection depth}\\[-5pt]
& \multicolumn{4}{c}{\hrulefill}&\multicolumn{4}{c@{}}{\hrulefill}\\
size &  Mean &  Min.  &  Max.  &  Variance &  Mean &  Min.  &  Max.  &  Variance \\
\hline
 25    & 0.3638 & 0.2571 & 0.4143 & 0.0018 & 0.0695 & 0.0286 & 0.1571 & 0.0017\\
 30    & 0.6571 & 0.5286 & 0.7429 & 0.0028 & 0.1333 & 0.0571 & 0.2429 & 0.0028\\
 35    & 0.7857 & 0.7000 & 0.8571 & 0.0017 & 0.1526 & 0.0429 & 0.2571 & 0.0044\\
 40    & 0.8390 & 0.7714 & 0.9000 & 0.0010 & 0.3010 & 0.1714 & 0.3714 & 0.0029\\
 45    & 0.8610 & 0.7714 & 0.9429 & 0.0019 & 0.4410 & 0.3571 & 0.5571 & 0.0041\\
 50    & 0.8990 & 0.8429 & 0.9286 & 0.0008 & 0.5124 & 0.4429 & 0.6000 & 0.0017\\
 55    & 0.9086 & 0.8714 & 0.9571 & 0.0006 & 0.5419 & 0.4429 & 0.6286 & 0.0021\\
\hline
\end{tabular*}
\vspace*{-6pt}
\end{table}

\begin{table}[b]
\caption{The mean, deviation and variance of the recognition rates by
TPD and PD with the Iris set} \label{tab:Iris_PD_TPD}
\begin{tabular*}{\textwidth}{@{\extracolsep{\fill}}lllllllll@{}}
\hline
\multirow{2}{31pt}{Training} & \multicolumn{4}{l}{Tensor projection depth}&\multicolumn{4}{l@{}}{Projection depth}\\[-5pt]
& \multicolumn{4}{c}{\hrulefill}&\multicolumn{4}{c@{}}{\hrulefill}\\
size &  Mean &  Min.  &  Max.  &  Variance &  Mean &  Min.  &  Max.  &  Variance \\
\hline
 10    & 0.9698 & 0.8571 & 1.0000 &0.0016 & 0.9476 & 0.7619 & 1.0000 & 0.0041\\
 15    & 0.9889 & 0.9524 & 1.0000 &0.0004 & 0.9841 & 0.9524 & 1.0000 & 0.0005\\
 20    & 0.9952 & 0.9048 & 1.0000 &0.0004 & 0.9921 & 0.8571 & 1.0000 & 0.0008\\
 25    & 0.9984 & 0.9524 & 1.0000 &0.0001 & 0.9984 & 0.9524 & 1.0000 & 0.0001\\
 30    & 0.9984 & 0.9524 & 1.0000 &0.0001 & 0.9984 & 0.9524 & 1.0000 & 0.0001\\
 35    & 1.0000 & 1.0000 & 1.0000 &0.0000 & 1.0000 & 1.0000 & 1.0000 & 0.0000\\
 40    & 1.0000 & 1.0000 & 1.0000 &0.0000 & 1.0000 & 1.0000 & 1.0000 & 0.0000\\
\hline
\end{tabular*}
\end{table}

From Table \ref{tab:Iris_PD_TPD} we can see that there is no apparent
difference between the two results. Therefore, data from vector spaces
can be converted  into tensors and we can perform the depth procession
with TPD.

\section{Discussion and conclusion}
\label{sec:conclusion}
In this paper,  tensor projection depth is
proposed as an extension of the definition of depth to tensor spaces.
We show that, according to the criteria \cite{Zuo2000a}, TPD satisfies
all four desirable properties. TPD has the advantages of avoiding the
curse of dimensionality and keeping the natural structures of the data
sets invariant. For  sparse samples, we use tensor PCA to remove their
null space and compute the TPD in the subspace. The numerical results
show that TPD performs better than PD for data which are naturally in
tensor form.

Data sets which are naturally in  vector form can also be processed
using TPD, which  converts the data into tensor form. Although such
processing will actually change the structure of the data sets to some
extent, numerical results show that there are no apparent differences
in the outcome. For some $(\mu, \sigma)$, such tensor-based processing
can effectively decrease the computational complexity of PD caused by
the dimensionality.

\begin{appendix}\label{app}
\section*{Appendix: Proofs}\label{sec:appendix}

\begin{pf*}{Proof of Theorem \ref{thm:properties}}
\textit{Convexity}. We will show that the outlying function
(\ref{eq:projection gen core tensor}) is still convex. Let
$\mathbf{X}_1,
\mathbf{X}_2\in\mathcal{R}^{n_1}\otimes\mathcal{R}^{n_2}$ be two
arbitrary points, $0<\lambda<1$, and for the point
$\mathbf{X}_0\doteq(1-\lambda)\mathbf{X}_1+\lambda\mathbf{X}_2$, we
have
\begin{eqnarray*}
&&|\mathbf{u}^T\mathbf{X}_0\mathbf{v}-\mu(\mathbf{u}^T\mathcal{X}\mathbf{v})|\\
&&\quad=\bigl|(1-\lambda)\bigl(\mathbf{u}^T\mathbf{X}_1\mathbf{v}-\mu(\mathbf{u}^T\mathcal{X}\mathbf{v})\bigr)+
\lambda\bigl(\mathbf{u}^T\mathbf{X}_2\mathbf{v}-\mu(\mathbf{u}^T\mathcal{X}\mathbf{v})\bigr)\bigr|\\
&&\quad\leq(1-\lambda)|\mathbf{u}^T\mathbf{X}_1\mathbf{v}-\mu(\mathbf{u}^T\mathcal{X}\mathbf{v})|+
\lambda|\mathbf{u}^T\mathbf{X}_2\mathbf{v}-\mu(\mathbf{u}^T\mathcal{X}\mathbf{v})|
\end{eqnarray*}
and
\begin{eqnarray*}
O(\mathbf{X}_0,\mathcal{X}) &=&\sup_{\|\mathbf{u}\|=\|\mathbf{v}\|=1}
\frac{|\mathbf{u}^T\mathbf{X}_0\mathbf{v}-\mu(\mathbf{u}^T\mathcal{X}\mathbf{v})|}
{\sigma(\mathbf{u}^T\mathcal{X}\mathbf{v})}\\
&\leq&\sup_{\|\mathbf{u}\|=\|\mathbf{v}\|=1}
\frac{(1-\lambda)|\mathbf{u}^T\mathbf{X}_1\mathbf{v}-\mu(\mathbf{u}^T\mathcal{X}\mathbf{v})|+
\lambda|\mathbf{u}^T\mathbf{X}_2\mathbf{v}-\mu(\mathbf{u}^T\mathcal{X}\mathbf{v})|}
{\sigma(\mathbf{u}^T\mathcal{X}\mathbf{v})} \\
&=&(1-\lambda)O(\mathbf{X}_1,\mathcal{X})+\lambda
O(\mathbf{X}_2,\mathcal{X}).
\end{eqnarray*}
Thus,
\[
\mathit{TPD}(\mathbf{X}_0,\mathcal{X})\geq
(1-\lambda)\mathit{TPD}(\mathbf{X}_1,\mathcal{X})+\lambda
\mathit{TPD}(\mathbf{X}_2,\mathcal{X}).
\]

\textit{Symmetry}. This is straightforward.

\textit{Affine invariance}. Suppose that $\mathbf{A}_{n_1\times n_1}$
and $\mathbf{B}_{n_2\times n_2}$ are any two non-singular matrices. We
then have
\[
O(\mathbf{AXB},\mathbf{A}\mathcal{X}\mathbf{B})
=\sup_{\|\mathbf{u}\|=\|\mathbf{v}\|=1}
\frac{|\mathbf{u}^T\mathbf{AXB}\mathbf{v}-\mu(\mathbf{u}^T\mathbf{A}\mathcal{X}\mathbf{B}\mathbf{v})|}
{\sigma(\mathbf{u}^T\mathbf{A}\mathcal{X}\mathbf{B}\mathbf{v})}.
\]
For  fixed $\mathbf{X}$, suppose that
\[
(\mathbf{u}_0,\mathbf{v}_0)=\arg\sup_{\|\mathbf{u}\|=\|\mathbf{v}\|=1}O(\mathbf{X},\mathcal{X}).
\]
Thus, if we fix $\mathbf{v}=\mathbf{v}_0$ and let
\[
\mathbf{u}_1=\arg\sup_{\|\mathbf{u}\|=1}
\frac{|\mathbf{u}^T\mathbf{A}(\mathbf{X}\mathbf{v}_0)-\mu(\mathbf{u}^T\mathbf{A}(\mathcal{X}\mathbf{v}_0))|}
{\sigma(\mathbf{u}^T\mathbf{A}(\mathcal{X}\mathbf{v}_0))},
\]
then, according to Theorem 2.1 in \cite{Zuo2003},
\[
\sup_{\|\mathbf{u}\|=1}
\frac{|\mathbf{u}^T\mathbf{A}(\mathbf{X}\mathbf{v}_0)-\mu(\mathbf{u}^T\mathbf{A}(\mathcal{X}\mathbf{v}_0))|}
{\sigma(\mathbf{u}^T\mathbf{A}(\mathcal{X}\mathbf{v}_0))}=\sup_{\|\mathbf{u}\|=1}
\frac{|\mathbf{u}^T(\mathbf{X}\mathbf{v}_0)-\mu(\mathbf{u}^T(\mathcal{X}\mathbf{v}_0))|}
{\sigma(\mathbf{u}^T(\mathcal{X}\mathbf{v}_0))}.
\]
Thus, $\mathbf{u}_1\mathbf{A}=\lambda\mathbf{u}_0$, where $\lambda\in
\mathbb{R}$, and we have
\begin{eqnarray*}
\sup_{\|\mathbf{v}\|=1}
\frac{|\mathbf{u}_1^T\mathbf{A}(\mathbf{X}\mathbf{v})-\mu(\mathbf{u}_1^T\mathbf{A}(\mathcal{X}\mathbf{v}))|}
{\sigma(\mathbf{u}_1^T\mathbf{A}(\mathcal{X}\mathbf{v}))}&=&\sup_{\|\mathbf{v}\|=1}
\frac{|\lambda\mathbf{u}_0^T(\mathbf{X}\mathbf{v})-\mu(\lambda\mathbf{u}_0^T(\mathcal{X}\mathbf{v}_0))|}
{\sigma(\lambda\mathbf{u}_0^T(\mathcal{X}\mathbf{v}))}\\
&=&\sup_{\|\mathbf{v}\|=1}
\frac{|\mathbf{u}_0^T(\mathbf{X}\mathbf{v})-\mu(\mathbf{u}_0^T(\mathcal{X}\mathbf{v}_0))|}
{\sigma(\mathbf{u}_0^T(\mathcal{X}\mathbf{v}))}.
\end{eqnarray*}
Therefore,
\[
\mathbf{v}_1=\arg\sup_{\|\mathbf{v}\|=1}
\frac{|\mathbf{u}_1^T\mathbf{A}(\mathbf{X}\mathbf{v})-\mu(\mathbf{u}_1^T\mathbf{A}(\mathcal{X}\mathbf{v}))|}
{\sigma(\mathbf{u}_1^T\mathbf{A}(\mathcal{X}\mathbf{v}))}=\mathbf{v}_0.
\]
Similarly,
\begin{eqnarray*}
\sup_{\|\mathbf{v}\|=1}
\frac{|\mathbf{u}_1^T\mathbf{A}\mathbf{X}\mathbf{v}-\mu(\mathbf{u}_1^T\mathbf{A}\mathcal{X}\mathbf{v})|}
{\sigma(\mathbf{u}_1^T\mathbf{A}\mathcal{X}\mathbf{v})}&=&
\sup_{\|\mathbf{v}\|=1}
\frac{|\mathbf{u}_1^T\mathbf{A}\mathbf{X}\mathbf{B}\mathbf{v}-\mu(\mathbf{u}_1^T\mathbf{A}\mathcal{X}\mathbf{B}\mathbf{v})|}
{\sigma(\mathbf{u}_1^T\mathbf{A}\mathcal{X}\mathbf{B}\mathbf{v})}\\
&=&\sup_{\|\mathbf{u}\|=\|\mathbf{v}\|=1}
\frac{|\mathbf{u}^T\mathbf{A}\mathbf{X}\mathbf{B}\mathbf{v}-\mu(\mathbf{u}^T\mathbf{A}\mathcal{X}\mathbf{B}\mathbf{v})|}
{\sigma(\mathbf{u}^T\mathbf{A}\mathcal{X}\mathbf{B}\mathbf{v})}\\
&=&
\frac{|\mathbf{u}_0^T\mathbf{A}\mathbf{X}\mathbf{B}\mathbf{v}_0-\mu(\mathbf{u}_0^T\mathbf{A}\mathcal{X}\mathbf{B}\mathbf{v}_0)|}
{\sigma(\mathbf{u}_0^T\mathbf{A}\mathcal{X}\mathbf{B}\mathbf{v}_0)}.
\end{eqnarray*}
The result then follows.

\textit{Monotonicity relative to deepest point}. Suppose that
$\mathbf{X}_1, \mathbf{X}_2,
\mathbf{X}_c\in\mathcal{R}^{n_1}\otimes\mathcal{R}^{n_2}$,
$\mathbf{X}_c$ is the deepest tensor and
$\mathbf{X}_1=\lambda\mathbf{X}_2+(1-\lambda)\mathbf{X}_c$,
$\lambda\in[0,1]$.  Then, since
\[
O(\mathbf{X}_1,\mathcal{X})\leq(1-\lambda)O(\mathbf{X}_2,\mathcal{X})+\lambda
O(\mathbf{X}_c,\mathcal{X}),
\]
we have
\[
O(\mathbf{X}_1,\mathcal{X})-\lambda
O(\mathbf{X}_c,\mathcal{X})\leq(1-\lambda)O(\mathbf{X}_1,\mathcal{X})\leq(1-\lambda)O(\mathbf{X}_2,\mathcal{X}).
\]
Thus, $O(\mathbf{X}_1,\mathcal{X})\leq O(\mathbf{X}_2,\mathcal{X})$ and
$\mathit{TPD}(\mathbf{X}_1,\mathcal{X})\geq \mathit{TPD}(\mathbf{X}_2,\mathcal{X})$, that
is, the tensor projection depth decreases monotonically along any ray
emanating from the deepest point.

\textit{Maximality at center}. Suppose that $F$ is $\theta$-symmetric
about a unique point $\mathbf{X}_c \in R^{n_1} \times R^{n_2}$. Then,
for any pair of unit vectors $\mathbf{u},\mathbf{v}$, we have
$\mu(\mathbf{u}^T\mathcal{X}\mathbf{v})=\mathbf{u}^T\mathbf{X}_c\mathbf{v}$
and the result follows.

\textit{Vanishing at infinity}. This is straightforward.
\end{pf*}

\begin{pf*}{Proof of Theorem \ref{thm:algorithm convergence}}
Define
\[
f(\mathbf{u},\mathbf{v})=\sup_{\|\mathbf{u}\|=\|\mathbf{v}\|=1}
\frac{|\mathbf{u}^T\mathbf{X}\mathbf{v}-\mu(F_n(\mathbf{u,v}))|}
{\sigma(F_n(\mathbf{u,v}))}.
\]
Let $\mathbf{u}_0$ be the initial value. Fixing $\mathbf{u}_0$, we get
$\mathbf{v}_0$ by solving the optimizations (\ref{eq:algorithm step 2})
and (\ref{eq:algorithm step 3}).

Likewise, fixing $\mathbf{v}_0$, we get $\mathbf{u}_1$ by solving the
optimization problem (\ref{eq:algorithm step 3}). Thus, we have
\[
f(\mathbf{u}_0, \mathbf{v}_0)\leq f(\mathbf{u}_1, \mathbf{v}_0).
\]
Finally, we get
\[
f(\mathbf{u}_0, \mathbf{v}_0)\leq f(\mathbf{u}_1, \mathbf{v}_0) \leq
f(\mathbf{u}_1, \mathbf{v}_1)\leq f(\mathbf{u}_2, \mathbf{v}_1)\cdots.
\]
Since $f$ is bounded, it converges.
\end{pf*}
\end{appendix}

\section*{Acknowledgements}
This work was supported by Grants from the Natural Science Fund of
China (Nos~60975038 and 60974124). We thank the two referees for their
careful reading and useful comments.

\printhistory

\end{document}